\documentclass[epsfig,11pt]{llncs}
\usepackage{amssymb,amsmath,graphicx,color}
\usepackage{amsfonts,fancyvrb}

\newcommand{\remove}[1]{}
\usepackage{epsfig,latexsym,graphicx}

\begin{document}

\title{\bf{The Exact Null Distribution of the \\Generalized Hollander-Proschan Type Test for NBUE Alternatives}}
\author{Kinjal Basu \and M.Z. Anis}
\institute{Indian Statistical Institute\\
203, B.T. Road, Kolkata - 700 108, India\\
Email : bst0707@isical.ac.in, zafar@isical.ac.in}
\maketitle

\begin{abstract}
In this note we derive the exact null distribution for the test statistic $\gamma_j^*(F_n)$ proposed by Anis and Mitra (2011) for testing exponentiality against NBUE alternatives. As a special case, we obtain the exact null distribution for the test statistic $K^*$ proposed by Hollander and Proschan (1975). Selected critical values for different size are tabulated for these two statistics. Some remarks concerning the benefits of using the exact distribution are made
\\\\
{\bf{Keywords}} : Mean residual life; Order statistics; Scale invariant.
\\\\
{\bf{MSC Classification: Primary: }}62E15, {\bf{Secondary: }}62G30, 62N05.
\end{abstract}

\section{Introduction}
Though the assumption of exponentiality is widely used in the theory of reliability and life testing (primarily because of its mathematical simplicity), it is often of interest to check for possible departure from exponentiality in the data. The notions of failure rate and the mean residual life are used to define a hierarchy of ageing classes which capture the effect of ageing, in some probabilistic sense, on the residual life. We have the well known classes of Increasing Failure Rate (IFR), Increasing Failure Rate on the Average (IFRA), New Better than Used (NBU), New Better than Used in Expection (NBUE) and Decreasing Mean Residual Life (DMRL), and their corresponding duals. We refer the interested reader to Barlow and Proschan (1975) and Bryson and Siddiqui (1969) for definitions and properties of these classes.

As an ageing class, the NBUE class is quite large. It is important in the study of replacement policies as shown by Marshall and Proschan (1972). Other important results for this class can be found in Barlow and Proschan (1975). Testing exponentiality against NBUE altenatives has received considerable attention during the last four decades. The first test was proposed by Hollander and Proschan (1975). Recently Anis and Mitra (2011) have genralized the Hollander - Proschan approach to propose a family of tests for NBUE alternatives. To the best of our knowledge, this is the only procedure which provides a family of test statistics for this class of life distributions. Other approaches for testing exponentiality against NBUE alternatives are compared in Anis and Basu (2011).

The purpose of this paper is to derive the exact null distribution of the generalized Hollander-Proschan type test for NBUE alternatives proposed by Anis and Mitra (2011). In fact we also obtain the exact null distribution for the test statistic $K^*$ proposed by Hollander and Proschan (1975) as a special case. These results are presented in Section 2. We present two tables giving the critical values for sample sizes $n = 2\;(1)\;25\;(5)\;100$ for different significance levels in Section 3. Finally, in Section 4 we discuss the results of a simulation study and conclude the paper.

\section{The Exact Null Distribution}
Let $\mathcal{E}$ be the class of exponential distributions with distribution function $F(x) = 1 - e^{-\lambda x}, x \geqslant 0$, where $\lambda$ is any positive number, typically known. Anis and Mitra (2011) consider the problem of testing 
\[\begin{split}
&H_0 : F \in \mathcal{E}\\
vs\;\; &H_1 : F \in NBUE - \mathcal{E}\\
\end{split}
\]
based on a random sample $X_1, X_2,\ldots, X_n$ of size $n$ from an absolutely continuous distribution $F$ with density $f$ and survival function $\bar{F}$. They measured the deviation of an NBUE distribution from exponentiality by the parameter
\[\gamma_j(F) = \int_0^\infty[\bar{F}(t)]^j[e_F(0) - e_F(t)]dF(t)
\]
where $j$ is a positive real number and $e_F(t) = \int_t^\infty\bar{F}(u)du/\bar{F}(t)$. Observe that when $j = 1$, we get $K\equiv \gamma_1(F) = \int_0^\infty[\bar{F}(t)][e_F(0) - e_F(t)]dF(t)$, the parameter considered by Hollander and Proschan (1975). 

To make the test scale invariant, Anis and Mitra (2011) consider
\begin{equation}
\gamma_j^*(F_n) = \frac{\gamma_j(F_n)}{\bar{X}_n}
\end{equation}where 
\[\begin{split}
\bar{X}_n &= \frac{\sum_{i=1}^nX_i}{n}\;\;\textnormal{and}\\
\gamma_j(F_n) &= \sum_{k=1}^nX_{(k)}\Biggl[\frac{1}{j}\left\{\biggl(\frac{n - k + 1}{n} \biggr)^{j+1} - \biggl(\frac{n - k}{n}\biggr)^{j+1}\right\} - \frac{1}{j(j+1)n}\Biggr]
\end{split}
\]
where $X_{(1)}\leqslant X_{(2)} \leqslant \ldots \leqslant X_{(n)}$ denotes the order statistics based on the random sample $X_1, \ldots, X_n$. For the special case $j = 1$, we get the Hollander - Proschan statistic $K \equiv \gamma_1(F_n) = \frac{1}{n^2}\sum_{k=1}^nX_{(k)}\{\frac{3n}{2} - 2k + 1\}$. Hollander and Proschan (1975) used $K^* = K/\bar{X}_n$ to make the test statistic scale invariant. Note that there is an error in the expression for $K$ as given by Hollander and Proschan (1975) as pointed out by Anis and Mitra (2011).

Both Anis and Mitra (2011) and Hollander and Proschan (1975) have shown that the asymptotic null distribution of their statistics is standard normal. Anis and Mitra (2011) made no mention regarding small sample. On the other hand, Hollander and Proschan (1975) observed that their statistic $K^*$ is just a linear function of the total time on the test statistic, and hence did not feel the need to furnish new small-sample critical points of $K^*$ because Barlow (1968) has (in his Table 3) given the percentile points of the total time on test statistic $\sum_{i=1}^{n-1}\{\tau(X_{(i)})/\tau(X_{(n)})\} = nK^* + \frac{1}{2}(n - 1)$; where $\tau(X_{(i)})$ is the total time on test to $X_{(i)}$ and is defined as $\tau(X_{(i)}) = \sum_{j=1}^iD_j$ where $D_j = (n - j + 1)(X_{(j)} - X_{(j-1)}), j = 1, \ldots, n-1$; with $X_{(0)} = 0$. 

Our aim here is to derive the exact null distribution of $\gamma_j^*(F_n)$ under the null hypothesis. By substituting $j=1$, we shall get the exact null distribution of the Hollander and Proschan (1975) test statistic $K^*$

To achieve this, we first write the statistic in terms of the normalized spacings, $ D_{k}\equiv ( n-k+1)(X_{(k)}- X_{(k-1)})$, in the following way:
\[ \gamma_{j}^{\star}\left( F_{n}\right)=\frac{\sum_{k=1}^{n}e_{k,n}^{\left( j\right) }D_{k}}{\sum_{k=1}^{n}D_{k}} \] 
where
\[e_{k,n}^{\left( j\right)}=\frac{1}{j}\left( \frac{ n-k+1}{n}\right)^{j} - \frac{1}{j(j+1)}\]
and assuming  $X_{(0)}=0$.

Since $\gamma_j^*(F_n)$ is scale invariant under $H_0$, we assume that $\lambda = \frac{1}{2}$, the exponentiality scale parameter. Then the $D_j$'s are independent $\chi _2^2$ variates. Now, by Theorem 2.4 of Box (1954),  we get, under the null hypothesis of exponentiality, the following result:

\begin{theorem}
Let X be a random variable with distribution function $ F_{X}\left( x\right) =1-exp\left( -\frac{x}{2}\right), x>0 $ and let $ X_{\left( 1\right) }\leq X_{\left( 2\right) }\leq\cdots\leq X_{\left( n\right) } $ denotes the order statistics based on the random sample  from $ X$. Then
\[P\left( \gamma_{j}^{\star}\left( X\right) \leq x\right)=1-   \sum_{i=1}^{n}\prod_{\begin {array}{l}
k=1\\
k\neq i\\
\end {array}}^{n} \left(\frac{e_{i,n}^{\left( j\right) }-x}{e_{i,n}^{\left( j\right) }-e_{k,n}^{\left( j\right) }}\right)I(x, e_{i,n}^{\left( j\right) })\]   
for $ e_{i, n}^{\left( j\right) }\neq e_{k, n}^{\left( j\right) }, $ for all $ i \neq k, $ for fixed $ n; $ where
\[ I\left( x, y\right) =\left\lbrace \begin{array}{cc}
1 & \mbox{~~~ if ~}  x < y  \\
0 & \mbox{~~~ if ~}  x \geq y.
\end{array} \right.\] 
\end{theorem}

Similarly, for the Hollander-Proschan statistic $K^*$, we have the following theorem:

\begin{theorem}
Let X be a random variable with distribution function $ F_{X}\left( x\right) =1-exp\left( -\frac{x}{2}\right), x>0 $ and let $ X_{\left( 1\right) }\leq X_{\left( 2\right) }\leq\cdots\leq X_{\left( n\right) } $ denotes the order statistics based on the random sample  from $ X$. Then
\[P\left( K^{\star} \leq x\right)=1-   \sum_{i=1}^{n}\prod_{\begin {array}{l}
k=1\\
k\neq i\\
\end {array}}^{n} \left(\frac{e_{i,n}-x}{e_{i,n}-e_{k,n}}\right)I(x, e_{i,n})\]   
for $ e_{i, n}\neq e_{k, n}, $ for all $ i \neq k, $ for fixed $ n; $ where
\[ I\left( x, y\right) =\left\lbrace \begin{array}{cc}
1 & \mbox{~~~ if ~}  x < y  \\
0 & \mbox{~~~ if ~}  x \geq y.
\end{array} \right.\]
and 
\[e_{i,n} = \frac{n + 2 - 2i}{2n}\]
\end{theorem}

\section{The Critical Values}
Anis and Mitra (2011) have suggested taking $j=0.25$. Accordingly, the critical values of $1.25\sqrt{1.5n}\gamma_j^*(F_n)$ are given in Table 1.
\begin{table}[h]
\footnotesize{
\caption{Critical Values for $1.25\sqrt{1.5n}\gamma_{0.25}^*(F_n)$}
\begin{center}
\begin{tabular}[c]{|c|c|c|c||c|c|c|}
\hline
n & $\alpha = 0.01$ & $\alpha = 0.05$ & $\alpha = 0.10$ & $\alpha = 0.90$ & $\alpha = 0.95$ & $\alpha = 0.99$\\
\hline
2&0.3679&0.4233&0.4919&1.5944&1.6635&1.7184\\
3&-0.2286&0.0152&0.1968&1.6125&1.7617&1.9614\\
4&-0.5940&-0.2080&0.0345&1.6098&1.7825&2.0605\\
5&-0.8368&-0.3566&-0.0848&1.6121&1.7992&2.1120\\
6&-1.0031&-0.4681&-0.1769&1.6090&1.8093&2.1475\\
7&-1.1286&-0.5557&-0.2519&1.6074&1.8181&2.1765\\
8&-1.2310&-0.6334&-0.3162&1.6019&1.8202&2.1964\\
9&-1.3191&-0.6893&-0.3657&1.5965&1.8240&2.2151\\
10&-1.3855&-0.7396&-0.4074&1.5924&1.8229&2.2271\\
11&-1.4484&-0.7912&-0.4522&1.5905&1.8286&2.2393\\
12&-1.4978&-0.8257&-0.4830&1.5859&1.8290&2.2522\\
13&-1.5396&-0.8601&-0.5154&1.5787&1.8258&2.255\\
14&-1.5835&-0.8933&-0.5413&1.5789&1.8300&2.2718\\
15&-1.6115&-0.9210&-0.5679&1.5742&1.8274&2.2754\\
16&-1.6497&-0.9470&-0.5887&1.5708&1.8268&2.279\\
17&-1.6764&-0.9737&-0.6129&1.5698&1.8280&2.2857\\
18&-1.6969&-0.9878&-0.6253&1.5656&1.8297&2.2923\\
19&-1.7223&-1.0073&-0.6436&1.5619&1.8270&2.2982\\
20&-1.7445&-1.0286&-0.6639&1.5570&1.8247&2.2939\\
21&-1.7699&-1.0435&-0.6762&1.5561&1.8266&2.3022\\
\hline
\end{tabular}
\end{center}
}
\end{table}
\begin{table}[h]
\footnotesize{
\begin{center}
\begin{tabular}[c]{|c|c|c|c||c|c|c|}
\hline
n & $\alpha = 0.01$ & $\alpha = 0.05$ & $\alpha = 0.10$ & $\alpha = 0.90$ & $\alpha = 0.95$ & $\alpha = 0.99$\\
\hline
22&-1.7829&-1.0599&-0.6908&1.5541&1.8250&2.3083\\
23&-1.8062&-1.0718&-0.7027&1.5481&1.8232&2.3106\\
24&-1.8198&-1.0857&-0.7150&1.5490&1.8239&2.3086\\
25&-1.8289&-1.0981&-0.7266&1.5432&1.8195&2.3094\\
30&-1.8922&-1.1544&-0.7770&1.5329&1.8162&2.3217\\
35&-1.9358&-1.1904&-0.8147&1.5286&1.8199&2.3338\\
40&-1.9729&-1.2264&-0.8448&1.5174&1.8120&2.3390\\
45&-1.9989&-1.2508&-0.8701&1.5068&1.8067&2.3444\\
50&-2.0171&-1.2702&-0.8909&1.5027&1.8045&2.3431\\
55&-2.0348&-1.2917&-0.9085&1.4935&1.7990&2.3469\\
60&-2.0573&-1.3094&-0.9271&1.4885&1.7931&2.3495\\
65&-2.0680&-1.3220&-0.9387&1.4832&1.7943&2.3555\\
70&-2.0788&-1.3363&-0.9541&1.4809&1.7933&2.3567\\
75&-2.1005&-1.3476&-0.9649&1.4777&1.7932&2.3599\\
80&-2.1035&-1.3567&-0.9739&1.4705&1.7855&2.3530\\
85&-2.1138&-1.3681&-0.9861&1.4690&1.7843&2.3593\\
90&-2.1171&-1.3760&-0.9938&1.4653&1.7843&2.3597\\
95&-2.1268&-1.3841&-0.9991&1.4619&1.7799&2.3575\\
100&-2.1342&-1.3947&-1.0106&1.4553&1.7753&2.3588\\
\hline
\end{tabular}
\end{center}
}
\end{table}
Recall that historically, the Hollander-Proschan statistic $K^*$ is the first test for exponentiality against NBUE alternatives. We give in Table 2, the critical values of $\sqrt{12n}\gamma_{1}^*(F_n)$. Note that here the exact distributions are used to obtain the critical values for $j = 2\;(1)\;25\;(5)\;60$, while for higher values of $n$, the critical values are based on $10^6$ simulated sample. It can be seen from the two tables that convergence to normality is extremely slow.
\begin{table}[h]
\footnotesize{
\caption{Critical Values for $\sqrt{12n}\gamma_{1}^*(F_n)$}
\begin{center}
\begin{tabular}[c]{|c|c|c|c||c|c|c|}
\hline
n & $\alpha = 0.01$ & $\alpha = 0.05$ & $\alpha = 0.10$ & $\alpha = 0.90$ & $\alpha = 0.95$ & $\alpha = 0.99$\\
\hline
2&0.0245&0.1227&0.2453&2.2037&2.3268&2.4251\\
3&-0.7183&-0.3687&-0.1067&2.1032&2.3667&2.7160\\
4&-1.0564&-0.5720&-0.2714&2.0028&2.3050&2.7862\\
5&-1.2379&-0.7018&-0.3916&1.9419&2.2529&2.7860\\
6&-1.3603&-0.7993&-0.4813&1.8915&2.2092&2.7745\\
7&-1.4492&-0.8711&-0.5456&1.8544&2.1788&2.7598\\
8&-1.5266&-0.9299&-0.5997&1.8209&2.1512&2.7439\\
9&-1.5792&-0.9765&-0.6391&1.7987&2.1320&2.7414\\
10&-1.6260&-1.0158&-0.6772&1.7712&2.1117&2.7255\\
11&-1.6627&-1.0478&-0.7061&1.7529&2.0922&2.7056\\
12&-1.7052&-1.0788&-0.7373&1.7354&2.0776&2.7052\\
13&-1.7329&-1.1002&-0.7569&1.7176&2.0608&2.6893\\
14&-1.7638&-1.1273&-0.7814&1.7028&2.0466&2.6804\\
15&-1.7796&-1.1428&-0.7946&1.6910&2.0398&2.6736\\
\hline
\end{tabular}
\end{center}
}
\end{table}

\begin{table}[h]
\footnotesize{
\begin{center}
\begin{tabular}[c]{|c|c|c|c||c|c|c|}
\hline
n & $\alpha = 0.01$ & $\alpha = 0.05$ & $\alpha = 0.10$ & $\alpha = 0.90$ & $\alpha = 0.95$ & $\alpha = 0.99$\\
\hline
16&-1.8005&-1.1593&-0.8133&1.6821&2.0296&2.6716\\
17&-1.8202&-1.1789&-0.8273&1.6686&2.0139&2.6524\\
18&-1.8332&-1.1920&-0.8414&1.6593&2.0075&2.6566\\
19&-1.8522&-1.2067&-0.8545&1.6496&1.9980&2.6526\\
20&-1.8669&-1.2176&-0.8659&1.6401&1.9915&2.6397\\
21&-1.8806&-1.2284&-0.8787&1.6292&1.9823&2.6296\\
22&-1.8933&-1.2432&-0.8892&1.6261&1.9773&2.6284\\
23&-1.8965&-1.2460&-0.8953&1.6190&1.9689&2.6161\\
24&-1.9134&-1.2598&-0.9049&1.6087&1.9637&2.6190\\
25&-1.9258&-1.2684&-0.9138&1.6051&1.9584&2.6052\\
30&-1.9580&-1.3014&-0.9460&1.5798&1.9355&2.5973\\
35&-1.9911&-1.3290&-0.9729&1.5581&1.9168&2.5840\\
40&-2.0147&-1.3521&-0.9944&1.5427&1.8994&2.5681\\
45&-2.0346&-1.3675&-1.0103&1.5252&1.8821&2.5506\\
50&-2.0507&-1.3840&-1.0239&1.5160&1.8751&2.5455\\
55&-2.0636&-1.3933&-1.0333&1.5046&1.8651&2.5332\\
60&-2.0841&-1.4104&-1.0500&1.4938&1.8511&2.5235\\
65&-2.0873&-1.4192&-1.0576&1.4884&1.8480&2.5153\\
70&-2.1029&-1.4283&-1.0668&1.4797&1.8419&2.5156\\
75&-2.1050&-1.4315&-1.0712&1.4703&1.8290&2.5040\\
80&-2.1167&-1.4416&-1.0810&1.4668&1.8275&2.5040\\
85&-2.1211&-1.4476&-1.0857&1.4639&1.8255&2.5003\\
90&-2.1321&-1.4539&-1.0938&1.4586&1.8193&2.4886\\
95&-2.1305&-1.4617&-1.0989&1.4556&1.8151&2.4877\\
100&-2.1407&-1.4615&-1.1022&1.4466&1.8088&2.4839\\
\hline
\end{tabular}
\end{center}
}
\end{table}

\section{Discussion}
Since Hollander and Proschan (1975) had suggested using the table given in Barlow (1968) for small sample sizes, we have attempted to see how good is this suggestion based on a simulation exercise. First note that, the critical values in Table 2, is derived by putting $j=1$ in $(1)$, by which the test statistic is a scaled version of $\gamma_1(F_n) = \frac{1}{n^2}\sum_{k=1}^nX_{(k)}\{\frac{3n}{2} - 2k + 1\}$. However, in their actual paper Hollander and Proschan (1975), have erroneously reported their test statistic to be a scaled version of $\frac{1}{n^2}\sum_{k=1}^nX_{(k)}\{\frac{3n}{2} - 2k + \frac{1}{2}\}$. 

We simulated a random sample of size $n$ from $Exp(1)$ distribution as the test is scale invariant. Based on this sample we calculated both the statistics $\gamma_1^*(F_n)$, and $K^*$. Note that both the test statistics are asymptotically same. We compare these two with the critical values given in Table 2 and also in Barlow (1968). We record whether the realization leads to an acceptance or rejection of the null hypothesis. Naturally, we would expect the test to accept the null hypothesis. We repeat this for $N=10^5$ times and see the proportion of times the test takes the correct decision. Thus we get the empirical size of the tests. This is reported in Table 3. All the simulations were done in PC platform using $\textnormal{MATLAB}^{\textregistered} 7.8.0$.

The nominal significane level is $5\%$. It is seen from Table 3, that the empirical size of the test based on $j=1$ is around $5\%$ while using the critical values in Table 2. Also note that the empirical size of the test based on $K^*$ is also around $5\%$ when we use the tables given in Barlow (1968). Note that the critical values given by Barlow (1968) for the total time on test statistic is just a linear function of the test statistic $K^*$ as proposed by Hollander and Proschan (1975). Thus they suggest to use these critical values. However, their test statistic $K^*$ contains an error as pointed out in Anis and Mitra (2011). Due to this error, the test statistic for $j=1$, i.e $\gamma_1^*(F_n)$ does not match $K^*$; though both are asymptotically same. In fact, if there was no mistake then, the critical values as suggested by Hollander and Proschan (1975) should have matched the critical values given in Table 2.

Thus here we see that because of the mistake committed in Hollander and Proschan (1975) but still using the critical values which are derived (independently) for a linear function of $K^*$, the test statistic, we get such good results. Hence, we recommend, that instead of using the errorneous test statistic for smaller values, it will be more appropriate to use the test given by $j=1$, namely, $\gamma_1^*(F_n)$ and the critical values given in Table 2. Furthermore note that if we compare the test statistic $\gamma_1^*(F_n)$ with the critical values given in Barlow (1968) we get dismal results (Column 3 of Table 3). Therfore, it is recommended more that we use the test statistic $\gamma_1^*(F_n)$ with the critical values in Table 2 for testing exponentiality versus NBUE family.

\begin{table}[h]
\footnotesize{
\caption{Size comparison of the Tests}
\begin{center}
\begin{tabular}[c]{|c|c|c|c|}
\hline
&\multicolumn {2}{|c|}{$\sqrt{12n}\gamma_1^*(F_n)$} & $\sqrt{12n}K^*$ \\
\hline
n & Using Table 2 & Using Barlow& Using Barlow\\
\hline
2&5.25&55.00&5.24\\
3&5.23&33.52&5.29\\
4&4.96&27.09&4.98\\
5&4.99&22.66&5.09\\
6&4.94&19.82&4.82\\
7&5.01&17.73&5.22\\
8&4.90&16.40&5.00\\
9&4.85&14.70&4.83\\
10&4.83&14.90&4.88\\
\hline
\end{tabular}
\end{center}
}
\end{table}


\begin{thebibliography}{99}
\bibitem{Anis}
Anis, M. Z. and Mitra, M (2011) A generalized Hollander-Proschen type test of exponentiality against NBUE alternatives. \textit{Statistics and Probability Letters}. \textbf{81};1;126-132
\bibitem{AnisBasu}
Anis, M. Z. and Basu, K. (2011). Monte Carlo comparison of test of exponentiality against NBUE alternative. In preparation.
\bibitem{Barlow}
Barlow, R. E. (1968). Likelihood ratio tests for restricted families of probability distributions. \textit{Annals of Mathematical Statistics.} \textbf{39}; 547-560
\bibitem{BP}
Barlow, R. E. and Proschan, F. (1975). \textit{Statistical Theory of Reliability and Life Testing}. Holt, Rinehart and Winston, New York.
\bibitem{Box}
Box, G. E. P. (1954). Some theorems on quadratic forms applied in the study of analysis of variance problems, I. Effect of inequality of variance in the one-way classification.   \textit{Annals of Mathematical Statistics.}
\textbf{25}, 290-302.
\bibitem{BSidd}
Bryson, M.C. and Siddiqui, M.M. (1969). Some criteria for aging. \textit{Journal of the American Statistical Association} \textbf{64}; 1472-14838
\bibitem{HP}
Hollander, M. and Proschan, F. (1975). Tests for the mean residual life. \textit{Biometrika}. \textbf{62} ; 585-593.
\bibitem{MP}
Marshall, A. W. and Proschan, F. (1972). Classes of distributions applicable in replacement with renewal theory implications. In: \textit{Proceedings of the 6th Berkeley  Symposium on Mathematical Statistics.} Eds: L.M. LeCam; J. Neyman and E. L. Scott. Vol. I. 395-415.
\end{thebibliography}
\end{document}